\documentclass[letterpaper, 10 pt, conference]{ieeeconf}
\usepackage{graphicx}

\usepackage{amsmath}
\usepackage{amsfonts}

\newtheorem{theorem}{Theorem}

\newtheorem{remark}{Remark}

\begin{document}

\title{\LARGE \bf Transverse Contraction Criteria for
  Existence,\\ Stability, and Robustness of a Limit Cycle\thanks{This work was supported by the Australian Research Council and the Boeing corporation.}}

\IEEEoverridecommandlockouts
\author{
Ian R. Manchester$^1$ \ \ Jean-Jacques E. Slotine$^2$ \\ \ \\
1: ACFR, School of Aerospace, Mechanical and Mechatronic Engineering, University of Sydney, Australia\\
2: Nonlinear Systems Laboratory, Massachusetts Institute of Technology, USA\\
ian.manchester@sydney.edu.au \ jjs@mit.edu}
\maketitle

\begin{abstract}
This paper derives a differential contraction condition for
the existence of an orbitally-stable limit cycle in an autonomous system.
This transverse contraction condition can be represented as a pointwise
linear matrix inequality (LMI), thus allowing convex optimization tools such as sum-of-squares programming to be
used to search for certificates of the existence of a stable limit cycle. Many
desirable properties of contracting dynamics are extended to this context, including
preservation of contraction under a broad class of
interconnections. In addition, by introducing the concepts of differential
dissipativity and transverse differential dissipativity,  contraction and transverse contraction can be established for  large scale systems via LMI conditions on component subsystems.
\end{abstract}

\section{Introduction}

Dynamic systems with periodic solutions are important in many areas of engineering, including biologically-inspired robot locomotion, phase-locked loops, vortex shedding from aircraft wings,  and combustion oscillations, to name just a few. In biology, oscillating systems seem to be the rule rather than the exception \cite{Rapp87}.

The basic question we address in this paper is the following: when does an autonomous system of the form
\begin{equation}\label{eq:basic_sys}
\dot x= f(x)
\end{equation}
have the property that all solutions starting from a particular set $K$ converge asymptotically to a unique limit cycle? It is well known that periodic solutions of an autonomous differential equation can never be asymptotically stable. This is clear from the fact that solutions which have initial conditions on the periodic orbit but offset in time will never converge. 

There is a long and distinguished history of research into limit cycles for nonlinear
systems. For example, the famous
result of Poincar\'e-Bendixson gives a very simple condition for
planar systems. An important generalization to monotone cyclic feedback systems was published in \cite{mallet1990poincare}, however this depends on quite a special system structure and there are many application areas where it does not apply. 

There are also interesting properties of the ``global'' structure of regions of attraction to periodic orbits. It is known that the region of attraction is a continuous deformation of a torus: the cartesian product of an open unit disc of dimension $n-1$, with a scalar circle coordinate \cite{wilson1967structure}. These are often referred to as ``transversal'' and ``phase'' coordinates, respectively. In all cases except possibly with $n=5$ it is guaranteed that the deformation is differentiable \cite{byrnes2010topological}, due to the recent resolution of the Poincar\'e Hypothesis by Perelman. Birkhoff gave necessary conditions for periodic solutions in terms of the existence of particular ``phase variables'', or associated differential one-forms \cite{birkhoff1927dynamical, byrnes2010topological}.

However, all of these conditions imply the existence of {\em at least one} limit cycle, but give no insight into the number of limit cycles, or their stability. In recent years many efficient computational methods for proving stability of equilibria of nonlinear systems have been proposed, using optimization methods to search for ``stability certificates'' such as Lyapunov functions and barrier certificates \cite{parrilo2003semidefinite}, \cite{tan2008stability}, \cite{prajna2007framework}, \cite{chesi2010lmi}. In previous papers, the first author and others have extended this computational approach to limit cycles analysis using ``transverse dynamics'' and sum-of-squares programming \cite{Shiriaev2008, Manchester10,Manchester10a,manchester2010regions}, however this method is not applicable when the system dynamics are uncertain, since uncertainty will generally change the location of the limit cycle in state space.

An alternative to Lyapunov methods is to search for a contraction metric  \cite{Lohmiller98}, \cite{aylward2008stability}. For the purposes of robust stability analysis of equilibria, an important difference is that a Lyapunov function must generally be constructed about a known equilibrium, whereas a contraction metric implies the existence of a stable equilibrium indirectly. This is particularly useful if the  equilibrium point may change location depending on the unknown dynamics.

Historically, basic convergence results on contracting systems can be traced back to the 1949 results of Lewis in terms of Finsler metrics \cite{lewis1949metric}, and results of Hartman \cite{hartman1961stability} and Demidovich \cite{demidovich1962dissipativity}. To our knowledge, contraction to limit cycles was first investigated using an identity metric by Borg  \cite{Borg}, and later by Hartman and Olech \cite{hartman1962global}.

In this paper, we introduce {\em transverse contraction}, extending the results of \cite{Borg}, \cite{hartman1962global} by exploiting generalized metrics and  system combination properties as in \cite{Lohmiller98}. We also give a nonlinear change of variables that  converts tranverse contraction to a linear matrix inequality (LMI) without conservatism. In Section \ref{sec:properties} we show that transverse contraction is preserved under several forms of interconnection with contracting systems. In Section \ref{sec:diffdiss} we introduce {\em differential dissipativity} and {\em transverse differential dissipativity}, as well as LMI conditions for each, giving a framework for optimization-based analysis of complex interconnections of nonlinear systems. Finally, in Section \ref{sec:ex}, we illustrate the applicability of the results on the Moore Greitzer jet engine model and for the identification of live neuron dynamics.

\section{Problem Setup and Preliminaries}

We assume that $f:K\rightarrow\mathbb R^n$ in \eqref{eq:basic_sys} is smooth and $x \in \mathbb R^n$, and that a unique solution of \eqref{eq:basic_sys} exists. We refer to the Jacobian of $f$ as $A(x):=\frac{\partial f}{\partial x}$. A set $K$ is called {\em strictly forward invariant} under $f$ if any solution of \eqref{eq:basic_sys} starting with $x(0)$ in $K$ is in the interior of $K$ for all $t> 0$. A periodic solution $x^\star$ is one for which there exists some $T>0$ such that $x^\star(t) = x^\star(t+T)$ for all $t$. Equilibria are trivially periodic for every $T$, but for oscillatory solutions -- which are our main concern -- there is some minimal time $T$ such that the above holds and this is referred to as the {\em period}. The {\em orbit} of a periodic solution is the set $\mathcal X^\star:=\{x: x=x^\star(t) \textrm{ for some } t\}$. Note that while non-trivial periodic {\em solutions} cannot be asymptotically stable, their {\em orbits} can be, and in this case we say that the solution is {\em orbitally stable} (see, e.g., \cite{Hale80}). Define a {\em time reparametrization} $\tau(t)$ as a smooth function $\tau:[0, \infty) \rightarrow [0, \infty)$ such that $\tau(t)$ is monotonically increasing and $\tau(t)\rightarrow\infty$ as $t\rightarrow \infty$.

\section{Contraction Conditions for Limit Cycles}\label{sec:contraction}

In this section we introduce a {\em transverse contraction} condition for an autonomous dynamical system
$\dot x = f(x), \ x \in M$, where $M$ is a smooth, compact $n$-dimensional manifold. The condition is given in terms of a function $V(x,\delta_x)$, where $x\in M$ and $\delta_x \in \mathbb R^n$, which induces a distance function similar to a Riemannian or Finsler metric \cite{bao2000introduction}.

For most of this paper, we will assume a Riemannian-like contraction metric $V(x,\delta_x) : = \sqrt{\delta_x'M(x)\delta_x}$ where $M(x)$ is positive-definite for all $x$, however the main results hold for more general structures such as Finsler metrics \cite{bao2000introduction, lewis1949metric, forni2012differential}. The following two theorems provide a generalization of the results of \cite{Borg}, \cite{hartman1962global}, which considered the case $V(x, \delta_x) = |\delta_x|^2$.

\begin{theorem}\label{thm:contr_zhuk}
Let $K\subset \mathbb R^n$ be compact, smoothly path-connected, and strictly forward invariant. If there exists a Finsler function $V(x,\delta)$ satisfying
\begin{equation}
\frac{\partial V}{\partial x}f(x) + \frac{\partial V}{\partial \delta_x}A(x)\delta_x\le -\lambda V(x,\delta_x),\label{eq:contr1}
\end{equation}
for all $\delta_x\ne 0$ such that $\frac{\partial V}{\partial \delta_x}f(x)=0$, then for every two solutions $x_1$ and $x_2$ with initial conditions in $K$ there exists time reparametrizations $\tau(t)$ such that $x_1(t)\rightarrow x_2(\tau(t))$ as $t\rightarrow\infty$. 
\end{theorem}

\begin{figure}
\begin{center}
\includegraphics[width=0.6\columnwidth]{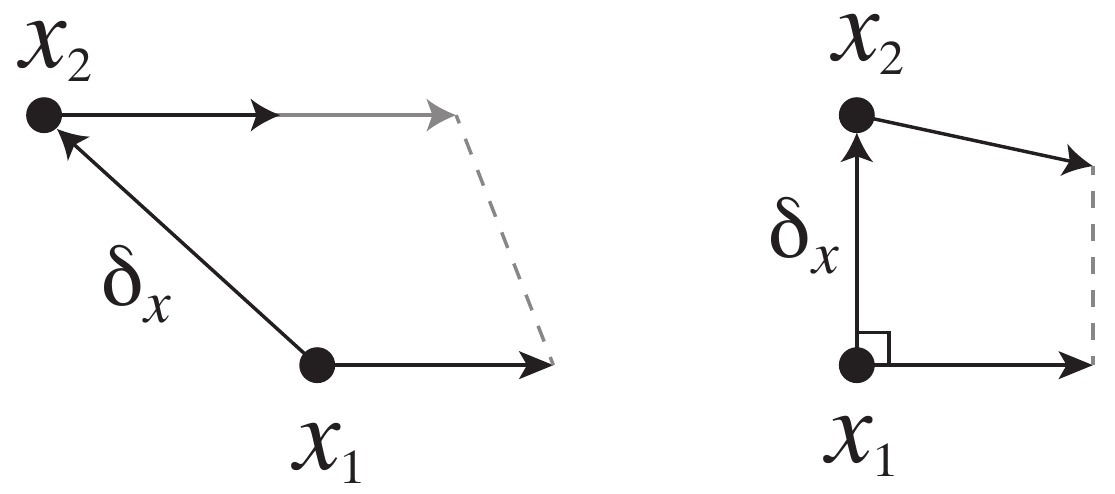}
\caption{On the left, an ``infinitesimal'' line segment joining $x_1$ and $x_2$ can be made to shrink be ``speeding up'' $x_2$ along its solution. On the right, the line segment is orthogonal to the derivative, so the system must be strongly contracting for the line segment to contract.}
\label{fig:transcon}
\end{center}
\end{figure}

\textbf{Proof:} The basic idea of the proof is illustrated in Figure \ref{fig:transcon}. Since the set $K$ is smooth and path-connected by definition, there
exists a smooth path between any two points $x_1\in K$ and $x_2\in K$ that remains in $K$. Such a
path can be considered as a smooth mapping $\gamma:[0,1]\rightarrow
\mathbb R^n$ with $\gamma(0) = x_1$ and $\gamma(1) = x_2$. We assume that paths are parametrized so that $\frac{\partial \gamma(s)}{\partial s} \ne 0$ for all $s$.

Denote by $\Gamma(x_1, x_2)$ the set of all such smooth paths between
$x_1$ and $x_2$ remaining in $K$ and associate with each a length
\[
L(\gamma) = \int_0^1 V\left(\gamma(s),\frac{\partial}{\partial s}\gamma(s)\right)ds
\]
and introduce the following Riemannian/Finsler-like distance between $x_1$ and $x_2$:
\begin{equation}\label{eq:distance}
d(x_1,x_2) = \inf_{\gamma \in \Gamma(x_1, x_2)} L(\gamma)
\end{equation}

The proof follows by showing that the distance $d(x_1, x_2)$ can by made to decrease by choice of time reparametrization, i.e. by speeding up or slowing down individual solutions along their phase portraits.

To this end, let us consider a path parametrized both in $s$ and time $t$: $\gamma(s,t)$, with the property that $\gamma(s,t_0)$ is the infimum in \eqref{eq:distance} for two points $x_1(t_0)$ and $x_2(t_0)$.
Now, let us introduce at every point $s\in[0, 1]$ and $t\ge t_0$ a ``speed scale'' $\alpha(s,t)>0$, which is assumed to be smooth in each argument. That is, at each point $\gamma(s, t)$ we have 
\[
\frac{d}{dt}\gamma(s,t) = \alpha(s,t)f(\gamma(s,t))
\]
with $\dot\tau(t)=\alpha(1,t)$ and $\alpha(0,t) = 1$

Now, by definition of the distance,
\[
\frac{d}{dt}d(x_1(t),x_2(\tau(t)) \le \int_0^1 \left[\frac{d}{dt}V\left(\gamma(s,t),\frac{\partial}{\partial s}\gamma(s,t)\right)\right]ds.
\]
Let us now consider, pointwise, the integrand in the right hand side of the above inequality.
\[
\frac{d}{dt}V\left(\gamma(s,t),\frac{\partial}{\partial s}\gamma(s,t)\right) = \frac{\partial V(x,\delta)}{\partial x}\dot x +\frac{\partial V(x,\delta_x)}{\partial \delta_x}\dot \delta_x 
\]
evaluated at $x = \gamma(s,t)$ and $\delta_x = \frac{\partial}{\partial s}\gamma(s,t)$, i.e. with
\begin{eqnarray}
\dot x &=& \alpha(s,t)f(\gamma(s,t)),\notag \\
\dot \delta_x &=& \frac{d}{dt}\frac{\partial}{\partial s}\gamma(s,t)
=\frac{\partial}{\partial s}(\alpha(s,t)f(\gamma(s,t)))\notag\\
&=&\frac{\partial \alpha}{\partial s}f(\gamma(s,t)))+\alpha(s,t)A(x)\frac{\partial \gamma}{\partial s}.\notag
\end{eqnarray}
Contraction under possible time-reparametrization follows from $\frac{d}{dt}V\left(\gamma(s,t),\frac{\partial}{\partial s}\gamma(s,t)\right)<0$ for all $s$. For this to hold for paths between all pairs of points, it is necessary that
\begin{equation}\label{eq:transcon_z}
\frac{d}{dt}V(x, \delta_x) = \frac{\partial V}{\partial x}f(x) +\frac{\partial V}{\partial \delta_x}\left(zf(x)+A(x)\delta_x\right)<0
\end{equation}
where the above has been normalized by $\alpha(s,t) > 0$ (which doesn't affect the sign) and 
where $z=\frac{1}{\alpha(s,t)}\frac{\partial \alpha}{\partial s}$ is a scalar. 

Since the time reparametrization is not specified, one interpretation is that $z$ is a ``control input'' which can be used to make the above inequality hold. Since it is affine in $z$, there are obviously ample choices of $z$ to satisfy this inequality as long as $\frac{\partial V}{\partial \delta_x}f(x)\ne 0$. The transverse contraction condition is simply that whenever $\frac{\partial V}{\partial \delta_x}f(x)=0$, \eqref{eq:transcon_z} is satisfied.\hfill$\Box$

\begin{remark}Stability under time reparametrization is sometimes referred to as {\em Zhukovsky stability} and has been used in several recent papers on limit
cycle stability, see e.g. \cite{yang2001remarks, Leonov06, Manchester10a} and apparently goes back to Poincar\'e in its essential argument \cite{Hale80}. 
It is known that systems satisfying such a property have limit cycles \cite{yang2001remarks}, but with the framework of contraction the proofs are simpler, so we give a proof here.
\end{remark}

\begin{theorem}\label{thm:main} If the conditions of Theorem 1 are satisfied, then all solutions starting with $x(0)\in K$ converge to a unique limit cycle.
\end{theorem}

\textbf{Proof:} since $K$ is
invariant and compact, it follows that $\Omega(x)$ exists and is a
compact subset of $K$. Furthermore, a clear implication of Theorem \ref{thm:contr_zhuk} is that all points in $K$
have the same $\omega$-limit set, which we denote $\Omega(K)$.

Pick a point $x^\star$ in $\Omega(K)$. By strict invariance, this is
an interior point of $K$. Assume $f(x^\star)\ne 0$, otherwise results of \cite{Lohmiller98} prove convergence to an equilibrium. Then one can
construct the hyperplane orthogonal to
$f(x^\star)$, which we denote $S$. We will prove convergence to a
limit cycle by
constructing a Poincar\'e map on $S$.

Since $f(\cdot)$ is smooth, for $x$ in some neighborhood $B$ of $x^\star$ we have
that $f(x)'f(x^\star)>0$, so in $B_S:=B\cap S$ solution curves are
transversal to $S$ and pass through it in the same direction as at
$x^\star$.

Since $x^\star$ is in the $\omega$-limit set for all points in $K$,
and $B_S$ is transversal,
the evolution of the system from any point $x(t)\in B_S$ eventually passes
through $B_S$ again, i.e. $x(t+s)\in B_S$ where $s>0$ depends on
$x$. This evolution can be represented by a Poincar\'e map
$T:B_S\rightarrow B_S$.

Take the distance between two points $d(x_1, x_2)$ on $B_S$ to be Riemannian metric
distance from Theorem \ref{thm:contr_zhuk}. Note that although the two points lie on the
$n-1$ dimensional set $B_S$, the curves joining them in the definition
of $d$ may pass out of the plane and through $n$-dimensional space. By
Thoerem \ref{thm:contr_zhuk}, we have that $d(T(x_1), T(x_2))<d(x_1, x_2)$. Hence $T$ is a
contractive map from $B_S$ onto itself, and by the Banach fixed point
theorem has a unique stable fixed point, which is its only limit
point so must be $x^\star$. By standard results on Poincar\'e maps
this implies that $x^\star$ is a point on a limit cycle, to which by
 all solutions converge, by Theorem \ref{thm:contr_zhuk}.\hfill$\Box$

\begin{remark} It can in fact be shown that convergence of $x_1(t)$ to the orbit of $x_2$ is {\em exponential} with rate $\lambda$, and that $\tau$ can be chosen to satisfy $\dot\tau(t) \rightarrow 1$ as $t\rightarrow\infty$, i.e. the system has {\em asymptotic phase}. We omit the details due to space restrictions.
\end{remark}

\begin{remark} Note that transverse contraction is a strictly weaker condition than contraction, so every contracting system is also transverse contracting. Hence the periodic solution to which a transverse contracting system converges may be {\em trivially} periodic, i.e. an equilibrium.
\end{remark}

\begin{remark} In \cite{wang2005partial}, \cite{Pham2007} and \cite{russo2011symmetries}, contraction transverse to a particular linear subspace was analyzed in the context synchronization. In this paper, contraction transverse to the system's vector field ensures asymptotically a form of ``synchronization'': in a periodic solution there is a single scalar variable (phase) that predicts all other states of the system. This concept may also be generalized to study higher-dimensional limit sets and non-autonomous systems.
\end{remark}


\subsection{Convex Formulation via Linear Matrix Inequalities}\label{sec:convex}

For the remainder of the paper we consider transverse contraction with a metric of the form $V(x, \delta_x) = \delta_x'M(x)\delta_x$. It will be shown in the next subsection that this class of metrics is sufficiently rich for testing orbital stability.

\begin{theorem}
A system $\dot x = f(x)$ is transverse contracting with rate $\lambda$ on a set $K$ if and only if there exists a function $\rho(x)\ge 0$  and a symmetric positive-definite matrix function $W(x)$ such that
\begin{equation}\label{eq:LMI}
 W(x) A(x) '+A(x) W(x)-\dot W(x) +\lambda W(x) - \rho(x) Q(x) \le 0
\end{equation}
for all $x\in K$, where $Q(x) := f(x)f(x)'$.
\end{theorem}

Note that this condition is linear in the unknown functions $W(x)$ and $\rho(x)$, i.e. it consists of a linear matrix inequality at each point $x$.

\textbf{Proof:} The following condition guarantees transverse contraction:
\[
\delta_x'\left(A(x) 'M(x)+M(x) A(x)+\dot M(x) +\lambda M(x))\right)\delta_x
\]
for all $\delta$ satisfying $\delta_x'M(x)f(x)=0$. If we reformulate this in terms of the {\em gradient} of the metric with respect to $\delta$: $\eta=M(x)\delta_x$, i.e. $\delta_x = M^{-1}(x)\eta=:W(x)\eta$
then
\begin{align}
&\delta_x'\left(A(x) 'M(x)+M(x) A(x)+\dot M(x) +\lambda M(x)\right)\delta_x \notag\\
 &= \eta'\left(W(x) A(x) '+A(x) W(x)-\dot W(x) +\lambda W(x))\right)\eta\notag
 \end{align}
since $\dot W(x) = \frac{d}{dt}(M^{-1}(x)) = -M^{-1}(x)\dot M(x) M^{-1}(x)$. Furthermore, the transversality condition $\delta_x'Mf=0$ is replaced by $\eta'f(x)=0$.

Now define matrix function $Q(x) := f(x)f(x)'$ which is rank-one and positive-semidefinite. This implies that the sets $\{\eta:\eta'f(x)=0\}$,  $\{\eta:\eta'Q(x)\eta=0\}$, and $\{\eta:\eta'Q(x)\eta\le 0\}$ are the same.

Transverse contraction with rate $\lambda$ can then be defined as the existence of a positive-definite matrix function $W(x)>0$ such that the following implication holds:
\begin{align}
&\eta'Q(x)\eta\le 0 \Rightarrow \notag\\
&\eta'\left(W(x) A(x) '+A(x) W(x)-\dot W(x) +\lambda W(x)\right)\eta\le 0 \notag
\end{align}
By the S-Procedure losslessness theorem \cite{yakubovich1971s}, the above implication is true if and only if there exists an $\rho(x)\ge 0$ such that
\[
 W(x) A(x) '+A(x) W(x)-\dot W(x) +\lambda W(x) - \rho(x) Q(x) \le 0
\]
which is the statement of the theorem.\hfill$\Box$

The above condition is convex and exact for each particular $x$. Such conditions can be verified over {\em regions} of the state space using sum-of-squares programming and positivstellensatz arguments \cite{parrilo2003semidefinite}, see \cite{aylward2008stability} for an exposition of this approach for the case of strong contraction.

\subsection{Generalized Jacobian and Transverse Linearization}\label{sec:singular}

The concept of a {\em generalized Jacobian} was introduced in \cite{Lohmiller98} for analysing contracting systems. Consider a nonsingular change of differential coordinates $\delta_z = \Theta(x)\delta_x$, then the dynamics in the new coordinates are given by $\dot\delta_z = F(x)\delta_z$ where the generalized Jacobian $F(x) := \Theta(x)A(x)\Theta(x)^{-1}+\dot\Theta(x)\Theta(x)^{-1}$. If such a change of coordinates exists such that $F(x)+F(x)' \le  -\lambda I$ then the system is contracting with rate $\lambda$. Furthermore, $M(x) = \Theta(x)'\Theta(x)$ is a valid contraction metric. Note that it is often easier to construct $\Theta(x)$ than a ``global'' change of coordinates $x\rightarrow z$.

A system is transverse contracting if there exists a differential change of coordinates such that $\delta_z(F+F')\delta_z < 0$ for all $\delta_z$ satisfying $\delta_z'\Theta(x)f(x) = 0$, where the latter condition follows from $\dot z = \Theta(x)\dot x = \Theta(x)f(x)$.

\begin{theorem}\label{thm:specialM}
If a system $\dot x = f(x)$ has a unique limit cycle to which all solutions starting in $K$ converge orbitally, then there exists a transverse contraction metric of the form $V(x, \delta_x) = \sqrt{\delta_x M(x)\delta_x}$ satisfying
\[
\frac{\partial V}{\partial x}f(x) + \frac{\partial V}{\partial \delta_x}A(x)\delta_x\le 0
\]
for all $\delta_x$ with strict inequality for $\delta_x$ satisfying $\frac{\partial V}{\partial \delta_x}f(x) = 0$. 
The generalized Jacobian is of the form
\[
F = \begin{bmatrix}0 & \star \\ 0 & F_\perp\end{bmatrix}
\]
where $F_\perp+F_\perp'<0$ and $F+F'$ has eigenvalues $0=\lambda_{max}>\lambda_2 \ge \lambda_3 ... \ge \lambda_n$.
\end{theorem}

\textbf{Proof:}
Here we include only a sketch of the proof due to space restrictions, the details are similar to the constructions in \cite{hauser1994converse, Leonov06, Manchester10a}.
In some toroidal neighbourhood $B$ of the limit cycle, there exists a smooth change of coordinates $x\rightarrow (\tau, x_\perp)$ where $\tau$ is a scalar phase variable along the cycle, and $x_\perp$ is an $(n-1)$-dimensional moving coordinate system orthogonal to $f(x)$. The differential system in these coordinates has the form
\[
\frac{d}{dt}\begin{bmatrix}\delta_\tau\\\delta_\perp\end{bmatrix} = \begin{bmatrix}0 & \star \\ 0& A_\perp(x)\end{bmatrix} \begin{bmatrix}\delta_\tau\\\delta_\perp\end{bmatrix} 
\]
Moreover, if the limit cycle is orbitally stable, there exists a Lyapunov function for the transversal part
\[
A_\perp(x)'M_\perp(x) + M_\perp(x)A_\perp(x) +\dot M_\perp(x)<0.
\]
A full metric is given by $|\delta_\tau|^2+\delta_\perp M_\perp(x)\delta_\perp$, which clearly satisfies the transverse contraction condition in $B$.
Since a solution from any point $x\in K$ converges to the limit cycle, there is a finite time after which it enters $B$. About this trajectory, a change of coordinates and Lyapunov function can be constructed via the method in \cite{Leonov06} satisfying the transverse contraction condition everywhere.

The construction of the generalized Jacobian comes from taking
\[
\Theta(x) = \begin{bmatrix}1 & 0 \\0 & \Theta_\perp(x)\end{bmatrix}\bar\Theta(x)
\]
where $\bar\Theta(x)$ is the Jacobian of the transformation $x\rightarrow (\tau, x_\perp)$ and $\Theta_\perp(x)$ satisfies $M_\perp(x) = \Theta_\perp(x)'\Theta_\perp(x)$.\hfill$\Box$

In the above, $A_\perp(x)$ is the transverse linearization that was used to construct Lyapunov functions for limit cycles in \cite{hauser1994converse} and \cite{Manchester10a}. Note that in those works, it was necessary for the limit cycle to be known and fixed to prove convergence, whereas transverse contraction decouples the question of convergence from knowledge of a particular solution.

\section{Properties of Transverse Contracting Systems}\label{sec:properties}

In many applications in which exact models are unavailable or very complex, it is desirable to characterize parameter ranges or interconnection structures over which the qualitative behaviour of the system remains the same. Engineering motivations are well known, but robustness analysis has also become of interest recently in biology, including as a measure of model validity \cite{morohashi2002robustness}. E.g., in \cite{morohashi2002robustness} robustness of limit cycles is assessed by gridding over parameter ranges and simulating the nonlinear system until convergence can be ascertained. Gridding and simulation becomes very expensive computationally for systems with large state dimension or many parameters, so alternative methods are desirable.

Feedback interconnections of oscillating systems with contracting systems may be of interest in many applications, for example control of robot arms \cite{williamson1998neural} or locomotion.

\subsection{Hierarchical Compositions of Systems}

A relatively simple application of the above theorem is to consider the composition of a contracting system and a transverse-contracting system.
\[
\dot x_1 = f_1(x_1), \ \ 
\dot x_2 = f_2(x_1, x_2)
\]

\begin{theorem} Suppose for each fixed $x_1$, $f_2$ is transverse contracting with metric $M_2(x_2,x_1)$, i.e.
\[
\delta_2'(F_2'M_{22}+M_2F_{22}+\frac{\partial}{\partial x_2} M_2f(x_2) +\lambda_2 M_2)\delta_2 \le 0
\]
for all $\delta_2$ satisfying $\delta_2'M_2f_2=0$ and $f_1$ is  strongly contracting in the sense of \cite{Lohmiller98}, i.e. there exists $M_1(x_1)$ such that
\[
F_1'M_1+M_1F_1+\frac{\partial}{\partial x_1} M_1f(x_1) +\lambda_1 M_1 \le 0
\]
then the composed system is transverse contracting, and hence has a unique stable limit cycle.
\end{theorem}
\textbf{Proof}
we prove this theorem by constructing a metric which decomposes as
\[
\delta'M_{c}\delta:=\delta_1'M_1\delta_1+\alpha\delta_2'M_2\delta_2=0
\]
which will be shown to verify the existence of a unique stable limit cycle. Let $x = [x_1' \ x_2']'$ and $f(x) = [f_1(x_1)' \ f_2(x_1,x_2)']'$. Since the contraction conditions are homogeneous with respect to $\delta$, and $\delta=0$ is trivial, it is sufficient to consider the case where $|\delta|=1$.

First, we note that since $f_2(x)\ne 0$ is $K$ and $K$ is compact, there exists an $\epsilon>0$ such that $|f_2(x)|\ge\epsilon$ for all $x\in K$.

Second, since system 1 is contracting and $K$ is compact, $f_1\rightarrow 0$ uniformly \cite{Lohmiller98}. The transversality condition for the metric $M_c$ is
\begin{equation}\label{eqn:compos_surf}
\delta_1'M_1f_1+\delta_2'M_2f_2=0.
\end{equation}
However, since $f_2$ is bounded below and $f_1$ converges uniformly to zero, the normal vector to the surface defined by \eqref{eqn:compos_surf} converges to that defined by $\delta_2'M_2f_2=0$ and hence the compact sets of $\delta$ of norm one satisfying these conditions converge uniformly.

Now, $d/dt [\delta'M_{c}\delta] = \delta'H\delta$ where $H$ decomposes into blocks corresponding to $\delta_1$ and $\delta_2$ like so:
 \[H=
\begin{bmatrix}
\alpha(F_1'M_1+M_1F_1+\dot M_1) & F_{21}'M_2+M_1F_{21}\\
F_{21}'M_1+M_1F_{21} & F_{22}'M_2+M_2F_{22}+\dot M_2
\end{bmatrix}.
\]
Consider the fixed $x_1^\star$ to which the contracting system $\dot x_1 = f_1(x_1)$ converges, so that $f_1(x^\star=0)$. Now, since system 1 is contracting, the upper-left block is negative definite and then by the Schur complement it follows that the maximum value of $\delta'H\delta$ on the subspace satisfying $\delta'Mf=0$ can be made strictly negative by choosing $\alpha$ sufficiently large.

Due to continuity of $H$ and the convergence of the sets of $\delta$, this implies that from any initial conditions there exists a finite time after which $\delta'H\delta<0$ for all $\delta$ satisfying \eqref{eqn:compos_surf}, which implies the existence of a unique stable periodic orbit by Theorems \ref{thm:contr_zhuk} and \ref{thm:main}.\hfill$\Box$

The opposite composition, a transverse-contracting system driving a contracting system clearly converges to a periodic solution due to natural input-to-state stability properties of contracting systems. In a sense, the second system can be considered as being driven by a periodic input \cite{Lohmiller98}.

\subsection{Robustness to Parametric Variation}

Suppose the system dynamics depend on some {\em parameter vector} $\theta$, i.e.
\[
\dot x = f(x,\theta).
\]
When studying robustness of equilibria of such systems, a widely-used method is to search for a parameter-dependent Lyapunov function (see, e.g., \cite{tan2008stability}).

In the context of the present paper, we assume that a particular set $K$ is robustly forward invariant -- which can be verified using the methods of \cite{tan2008stability} -- then robust existence of a single globally stable (within $K$) limit cycle is ensured if one can find a parameter-dependent contraction metric $M(x,\theta)$ which satisfies
\[
\delta'(\dot M(x,\theta) + 2F(x,\theta)'M(x,\theta) + \lambda M(x,\theta) )\delta \le 0
\]
for all $\delta$ such that $\delta'M(x,\theta)f(x,\theta)=0$, and for all $x\in K$ and $\theta$ in some set $\Theta$, where $K$ is a forward invariant set.

Note that this condition can be expressed as a parameter-dependent LMI as in \eqref{eq:LMI}, and verified via either sum-of-squares \cite{parrilo2003semidefinite} or sample-based methods \cite{calafiore2006scenario}.

\subsection{Skew-Symmetric Feedback Interconnection}
In this section and the next one we consider feedback interconnections of two systems of the form:
\begin{equation}
\dot x_1 = f_1(x_1,x_2), \ \ 
\dot x_2 = f_2(x_1,x_2).\label{eq:skewsym}
\end{equation}

\begin{theorem}
Suppose System 1 is partially contracting with respect to $x_1$, i.e. there exists a differential change of coordinates $\Theta_1(x_1)$ such that $F_1:=\Theta_1\frac{\partial f_1}{\partial x_1}\Theta_1^{-1}+\dot\Theta_1\Theta_1^{-1}$ satisfies $F_1+F_1'<0$.

Suppose also that System 2 is partially transverse contracting with respect to $x_2$, i.e. by Theorem \ref{thm:specialM} there exists a differential change of coordinates $\Theta_2(x_2)$ such that $F_2:=\Theta_2\frac{\partial f_2}{\partial x_2}\Theta_2^{-1}+\dot\Theta_2\Theta_2^{-1}$ satisfies $F_2+F_2'\le 0$ and $\delta_2(F_2+F_2')\delta_2<0$ when $\delta\ne 0$ satisfies $\delta_2\Theta_2f_2=0$.

Define $G_{12}:=\Theta_1\frac{\partial f_1}{\partial x_2}\Theta_2^{-1}$ and $G_{21}:=\Theta_2\frac{\partial f_2}{\partial x_1}\Theta_1^{-1}$ and suppose $G_{12} = -kG_{21}'$ for some $k>0$, then the interconnection \eqref{eq:skewsym} is transverse contracting.
\end{theorem}

\textbf{Proof}: Let $f = [f_1' \ f_2']$ and $x = [x_1' \ x_2']$. We will make use of the differential change of coordinates
\[
\Theta = 
\begin{bmatrix}\Theta_1&0\\0&\sqrt{k}\Theta_2\end{bmatrix}
\]
and define $F:=\Theta\frac{\partial f}{\partial x}\Theta^{-1}+\dot\Theta\Theta^{-1}$. The interconnection is transverse contracting if $\delta'(F+F')\delta<0$ for all $\delta$ such that $\delta'\Theta f
=0$.

First, we decompose $\delta = [\delta_1' \ \delta_2']'$ matching the decomposition of $x$, after some simple algebra we see that the off-diagonal terms cancel, so the transverse contraction condition is
\begin{equation}\label{eq:skewsymmcontr}
\delta_1'(F_1+F_1')\delta_1 +\delta_2'(F_2+F_2')\delta_2<0,
\end{equation}
for all $\delta_1, \delta_2$ not both zero satisfying $\delta_1'\Theta_1f_1+\sqrt{k}\delta_2'\Theta_2f_2=0$. Let us consider two cases:

Case 1: $\delta_2 =0$.
In this case the transversality condition $\delta'\Theta f=0$ reduces to $\delta_1'\Theta_1f_1 = 0$ and contraction is $\delta_1'(F_1+F_1')\delta_1 <0$. So this reduces to the assumed transverse contraction of System 1.

Case 2: $\delta_2 \ne 0$
Condition \eqref{eq:skewsymmcontr} is satisfied because $\delta_2(F_2+F_2')\delta_2<0$ for nonzero $\delta_2$ and  $F_1+F_1$ is negative semidefinite, hence $\delta'(F+F')\delta<0$.\hfill $\Box$

\subsection{Bounded Feedback Interconnections}

A more general theorem was presented in \cite{wang2005partial} for contracting systems. Here we discuss how it extends to transverse contraction. Suppose we have a general feedback interconnection, and construct $F$ as above. Define

\[
F_s := F+F' = \begin{bmatrix}F_{1s} & G_s\\ G_s' & F_{2s}\end{bmatrix}\begin{bmatrix}\delta_1\\ \delta_2\end{bmatrix} =:F\delta
\]
where $F_{1s} := F_1+F_1'$ and $F_{2s} := F_2+F_2'$ and $G_s := \Theta_1\frac{\partial f_1}{\partial x_2}\Theta_2^{-1}+\left(\Theta_2\frac{\partial f_2}{\partial x_1}\Theta_1^{-1}\right)'$.

Suppose system 1 is transverse contracting, so $F_{1s}\le 0$ and $z'F_{1s}z<0$ for all $z'\Theta_1f_1=0$. In \cite{wang2005partial} the Schur complement was used to derive conditions for contraction:
\[
F_{1s}\le G_sF_{2s}^{-1}G_s' \Leftrightarrow F_s \le 0.
\]
Note that in the case of transverse contracting systems, $z'F_{1s}f_1(x)z$ when $z'\Theta f=0$  and $z'F_{2s}z<0$ otherwise. Since $F_{2s}$ is nonsingular,  for the inequality on the left hand side to hold, it must be the case that $G_s'\Theta f = 0$. 
A very simple condition for $G_s=0$, which is equivalent to the skew-symmetric condition in the previous section, i.e. $\Theta_1\frac{\partial f_1}{\partial x_2}\Theta_2^{-1}=-\left(\Theta_2\frac{\partial f_2}{\partial x_1}\Theta_1^{-1}\right)'$

Another sufficient condition for $Gf_1(x)=0$ would be for both of these terms to be zero. For the first term, this implies that perturbations in System 2 only affect the {\em transversal} states of System 1, not the phase. For the second term, this means that perturbations in the phase of system 1 do not affect system 2. This would correspond to a decomposition of System 1 into a phase and transversal system, only the latter of which interacts with System 2.

Suppose that $G'\Theta f_1=0$ then a sufficient condition for transverse contraction of the interconnection is
\[\lambda_2(F_{1s})\lambda_{max}(F_{2s})<\sigma^2(G)\]
by a similar argument to \cite{wang2005partial}.
Note that $\lambda_2(F_{1s})$ is the rate of transverse contraction of System 1 and $\lambda_{max}(F_{2s})$ is the exponential rate of contraction of System 2.

\subsection{Robustness to Bounded Disturbance}

Consider the global coordinates $x_\perp, \tau$ -- either implicitly or explicitly defined. Since $\tau\in S^1$ the dynamics of $x_\perp$ can be considered a periodic differential equation with a transformation of time. This makes it clear that any internal perturbation in $f$ which keeps $\dot\tau>0$ and $F_\perp(x)$ contracting still results in a limit cycle (c.f. above).

Bounded external perturbations will also have bounded effect on behavior. Denote $x^\star$ the periodic orbit of a transverse contracting system $\dot x = f(x)$. Letting $R(x) = \min_{\tau}\int_{x^\star(\tau)}^{x} V(\gamma(s), \frac{\partial \gamma}{\partial s})ds$ we have
\[
\dot R + \lambda R\le 0
\]
Consider a bounded external disturbance, i.e.
$\dot x = f(x)+d(t)$, where $|d|\le d_{\max}$, then we have
\[
\dot R+ \lambda R\le |\Theta d(t)|
\]
so after exponentially-forgotten transients, the perturbed system is within a ball of radius $R$ around the original limit cycle. For further details on such analysis, see \cite{Lohmiller98}.

\section{Differential Dissipativity and \\Transverse Differential Dissipativity}\label{sec:diffdiss}

Methods related to dissipation inequalities are central to quantitative results in systems analysis, including input-output methods such as small-gain and passivity \cite{desoer1975feedback}, robust control design \cite{petersen2000robust}, and integral quadratic constraints \cite{megretski1997system, IQCchapter}. In this section, we introduce concepts of differential dissipativity, closely related to incremental small gain and passivity \cite{desoer1975feedback}. 

Roughly speaking, a system is differentially dissipative if the {\em linearization} along every solution is dissipative, however the results are exact and global, not local. The concept has been used several times before -- though not under that name -- in constructing small gain theorems for contracting systems \cite{jouffroy2003simple} and in bounding the simulation error of identified models \cite{tobenkin2010convex, manchester2011}.

For this section we consider systems with external inputs and outputs:
\begin{equation}\label{eqn:inputsys}
\dot x = f(x,w), y = g(x,w)
\end{equation}
which has the differential system:
\begin{equation}\label{eqn:inputdiff}
\dot \delta_x = A(x)\delta_x+B(x)\delta_w, \   \delta_y = C(x)\delta_x+D(x)\delta_w, 
\end{equation}
where $A(x):=\frac{\partial f}{\partial x}, B(x) := \frac{\partial f}{\partial w}, C(x): = \frac{\partial g}{\partial x}, D(x) : = \frac{\partial g}{\partial w}$.

A statement about {\em differential dissipativity} relates the system \eqref{eqn:inputsys}, \eqref{eqn:inputdiff} to a particular form $\sigma(x,w,\delta_x, \delta_w)$ which in applications is usually quadratic in $\delta_x, \delta_w$. In particular, along all solutions of \eqref{eqn:inputsys}, the differential system \eqref{eqn:inputdiff} satisfies
\begin{equation}\label{eqn:diss_in}
\int_0^T\sigma(x,w,\delta_x, \delta_w)dt \ge -\kappa(x(0),\delta_x(0))
\end{equation}
for all $T>0$ and for some $\kappa:TM\rightarrow \mathbb R$. A shorthand notation for this is $\sigma(x,w,\delta_x, \delta_w) \succ 0$, c.f. the notion of a ``complete IQC'' in \cite{IQCchapter}.

For example, one can define differential versions of the classical small-gain condition with $\sigma_{\gamma}=\gamma |\delta_w|^2-|\delta_y|^2$ and passivity with $\sigma_p=\delta_w'\delta_y$, where the latter assumes the input and output have matching dimensions. 

Inspired by IQC analysis \cite{IQCchapter}, if a number of system properties are encoded in dissipativity relations of the form $\sigma_i \succ 0, i = 1, 2, ... p$, then a desired property (e.g. stability or bounded gain) encoded as $\sigma^\star \succ 0$, and then one searches for constants $\tau_i \ge 0, i = 1, 2, ... p$ satisfying $
\sigma^\star-\sum_{i=1}^p\tau_i\sigma_i\succ 0.
$

For system evolution on an invariant compact set, taking $\sigma^\star := -|\delta_x|^2\succ 0$ implies contraction, since it implies that $\delta_x$ converges to zero via Barbalat's lemma \cite{khalil2002nonlinear}. Differential contraction versions of the small-gain theorem and the passivity theorem are special cases of this formulation.

For a system of the form \eqref{eqn:inputsys}, a sufficient condition for  \eqref{eqn:diss_in} is the existence of a metric function $V(x,\delta_x) = \delta'M(x)\delta>0$ such that
\begin{equation}\label{eq:tdissineq}
\frac{d}{dt}V(x,\delta_x) \le \sigma(x,w,\delta_x, \delta_w)
\end{equation}
where the path integral of $V$ plays the role of an incremental storage function between solutions. 

We define a system as {\em transverse differentially dissipative} (TDD) with a supply rate $\sigma(x,w,\delta_x, \delta_w)$ if \eqref{eq:tdissineq} holds for all $\delta_x$ such that $\frac{\partial V}{\partial \delta_x}f(x,w)=0$.

We give the following theorem, which can easily be extended to more than two system.
\begin{theorem} Given two systems $\dot x_1 = f_1(x_1,w_1)$ and $\dot x_2 = f_2(x_2, w_2)$, and consider the interconnection $w_1 = g_2(x_2,w_2), w_2 = g_1(x_1, w_1)$. Suppose System 1 is transverse differentially dissipative with respect to supply rate $\sigma_1(x_1,w_1,\delta_{x1}, \delta_{w1})$ and System 2 satisfies $\sigma_1(x_1(t),w_1(t),\delta_{x1}(t), \delta_{w1}(t))\ge 0$ for all $t$. Then if there exists nonnegative constants  $\tau_1, \tau_2$ such that $0<\tau_1\sigma_1(x_1,w_1,\delta_{x1}, \delta_{w1})+\tau_2\sigma_1(x_1,w_1,\delta_{x1}, \delta_{w1})$ on a forward-invariant set of the interconnected system, then the interconnection is transverse contracting and has a unique stable periodic solution.
\end{theorem}
The proof of this theorem follows standard S-Procedure arguments in robust control theory \cite{petersen2000robust, IQCchapter}.

For example, for a dynamic system in feedback with a time-varying but non-dynamic mapping $w = \Delta(y,t)$ where $\Delta$ is slope-restricted with respect to $y$, one can choose $\sigma_1(x,w,\delta_x,\delta_w)  \sigma_c(\delta_y,\delta_w):= (\delta_y-\alpha \delta_w)(\beta\delta_w-\delta_y)$ and  $\sigma_2 = -\sigma_c(\delta_w,\delta_y)$. In doing so, we recover a differential form of the circle criterion that proves existence of a limit cycle in feedback with sector-bounded and slope-restricted nonlinearities.

\subsection{Linear Matrix Inequalities for DD and TDD}

The convex formulation from Section \ref{sec:convex} can be extended to differential dissipativity conditions where the supply rate has the form
\[
\sigma = \delta_x'H(x,u)\delta_x + 2\delta_x'N(x,u)\delta_u + \delta_u'R(x,u)\delta_u
\]
as long as $H(x,u)$ is negative semidefinite, which is the case for common supply rates such as passivity and small gain. Note that it is necessary that $R(x,u)$ to be positive semidefinite for a lower bound to exist in \eqref{eqn:diss_in}.

Using the S-Procedure, the following condition is equivalent to transverse differential dissipativity:
\begin{align}
 \eta'(-WA'-AW+\dot W+\rho Q+WHW)\eta & \notag\\+ 2\eta'(-B+WN)\delta_u+\delta_uR\delta_u &\le 0
\end{align}
where we have dropped dependence of matrices on $x$ and $u$ for the sake of space and clarity. 
Note that although this inequality is quadratic in $W$ it is still convex (when $H\le 0$) and it can be linearized via a Schur complement to give the following condition:
\[
\begin{bmatrix}
-WA'-AW+\dot W+\rho Q & -B+WN& W\\
 -B'+WN'&R& 0\\
 W & 0 &-H
 \end{bmatrix}\ge 0.
\]
Here the matrices $A, B, H, N, R, Q$ are specified by the system description, and the decision variables are the certificate functions $W, \rho$. The above matrix inequality is clearly linear in the decision variables, and is therefore amenable to search via convex optimization.

\section{Application Examples}\label{sec:ex}

\subsection{Moore Greitzer model of combustion oscillation}

The Moore-Greitzer model, a simplified model of surge-stall dynamics of a jet engine  \cite{moore1986theory}, has motivated substantial development in nonlinear control design (see, e.g.,  \cite{krstic1995nonlinear} and references therein). In \cite{aylward2008stability}, sum-of-squares programming was applied for automated construction of verification of contraction metrics. The following form of the Moore Greitzer model was examined, with $\delta$ considered an uncertain parameter:
\[
\begin{bmatrix}\dot\phi \\ \dot\psi \end{bmatrix}  = \begin{bmatrix} -\psi -\frac{3}{2}\phi^2 - \frac{1}{2}\phi^3+\delta\\ 3\phi - \psi \end{bmatrix}.
\]
Contraction, and hence existence of a stable equilibrium, was established that values of $\delta$ with $|\delta|<1.023$ using a contraction metric with each element a degree-six polynomial. In fact the system is also contracting for values of $\delta>1.023$, but at $\delta\approx -1.023$ a Hopf bifurcation occurs.

Using the S-procedure formulation for transverse contraction from Section \ref{sec:convex} of the present paper, we have established that for values of $\delta<-1.023$ the Moore Greitzer model exhibits stable oscillations. 

Let $H(x) = A(x)W(x)+W(x)A(x)'-\dot W(x)+\lambda W(x)$, and let $\Sigma[x]$ denote the set of sum-of-squares polynomials in $x$, and $\Sigma_n[x]$ denote the set of $n\times n$ matrices verified positive semidefinite via sum-of-squares i.e. matrices $R(x)$ satisfying $y'R(x)y\in\Sigma[x,y]$.

Using a positivstellensatz construction \cite{parrilo2003semidefinite} we derive the following conditions for transverse contraction, restricted to a set $K$ which is a disc of radius $\rho$ with a small region around the unstable equilibrium deleted.
\begin{eqnarray*}
W(x)  -(f(x)'f(x)-0.1)L_1(x)&&\\-(\rho-x'x)L_2(x) &\in& \Sigma_n[x],\\
-H(x) -\alpha(x)f(x)f(x)'&&\\-(f(x)'f(x)-\epsilon)L_3(x)-(\rho^2-x'x)L_4(x) &\in& \Sigma_n[x],\\
L_1(x), L_2(x), L_3(x), L_4(x), &\in& \Sigma_n[x],\\
\alpha(x) &\in& \Sigma[x].
\end{eqnarray*}
We found that these conditions could be verified with $\rho = 10, \epsilon = 0.1$, and $W(x)$ a matrix of degree-four polynomials, and $L_i(x), \alpha(x)$ degree-two. The MATLAB code used to verify these conditions has been made available online \cite{MGexample}.

\begin{figure}
\begin{center}
\includegraphics[width=0.8\columnwidth]{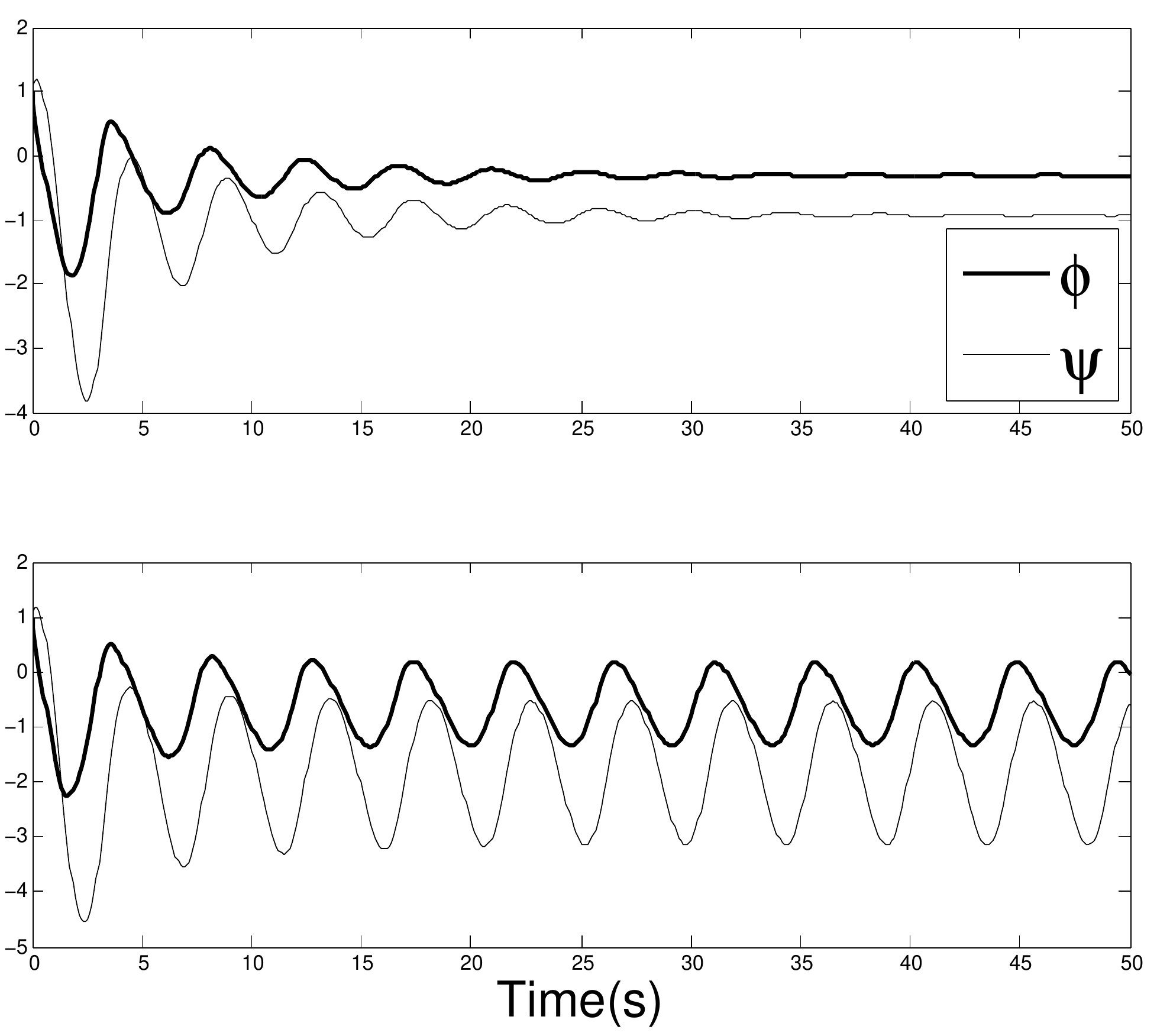}
\caption{Moore-Greitzer jet engine model response with $\delta = -0.8$ (left) and $\delta = -1.2$ (right).}
\label{fig:hist}
\end{center}
\end{figure}

\subsection{Identification of Oscillating Systems: Live Neurons}

Identifying nonlinear models with stable oscillations is a highly challenging problem. A new framework for nonlinear state-space system identification was introduced in \cite{tobenkin2010convex} which can be used to guarantee stability of identified models, and in \cite{manchester2011} this method was extended to allow stable limit cycles, although that paper did not contain strong theoretical claims. The problem is: given a measured set of data points $\tilde x,\dot {\tilde x}$ find a stable nonlinear differential equation that reproduces the data. The proposed method searches over a very flexible class of models: $
E(x)\dot x=f(x)
$
where $E(x)$ and $f(x)$ are matrices of polynomials, and $E(x)$ is nonsingular. A special form of a metric was proposed:
\[
M(x) = \Pi(x)'E(x)'QE(x)\Pi(x)
\]
where $Q$ is a positive definite matrix and $\Pi(x)$ is a projection on to the subspace orthogonal to $\dot x$. The main result of  \cite{manchester2011} is a reformulation of the problem of joint search for dynamics and metric -- i.e. $E(x), f(x)$ and $Q$ -- as a convex optimization problem (a sum-of-squares program). 

In \cite{manchester2011} this method was used to accurately identify dynamic models of live rat hippocampal neurons in culture, including both contracting sub-threshold dynamics and orbitally stable periodic ``spiking''.

The results on transverse contraction in the present paper lend theoretical justification to this procedure, showing that such a metric does in fact enforce the existence of stable limit cycles for the model, with some caveats due to approximations used in \cite{manchester2011}. A more complete discussion of this will follow in another publication.

\bibliographystyle{IEEEtran}
\bibliography{elib}

\end{document}